\documentclass[12pt]{amsart}

\usepackage[latin1]{inputenc}
\usepackage{subfigure}
\usepackage{amsmath}
\usepackage{color}
\usepackage{amsfonts, amssymb, hyperref, enumerate}
\usepackage{graphics}

\numberwithin{equation}{section}
\newcommand{\N}{\mathbb{N}}
\newcommand{\Z}{\mathbb{Z}}

\newcommand{\C}{\mathbb{C}}
\renewcommand{\S}{{\mathcal S}}

\newcommand{\beq}{\begin{equation*}}
\newcommand{\eeq}{\end{equation*}}
\newcommand{\be}{\begin{equation}}
\newcommand{\ee}{\end{equation}}

\theoremstyle{plain}
\newtheorem{theorem}{Theorem}[section]

\newtheorem{prop}[theorem]{Proposition}

\newtheorem{lemma}[theorem]{Lemma}
\newtheorem{cor}[theorem]{Corollary}

\newtheorem{example}[theorem]{Example}

\newtheorem{rem}[theorem]{Remark}
\newtheorem{rems}[theorem]{Remarks}

\newcommand{\al}{\alpha}
\newcommand{\la}{\lambda}

\newcommand{\lan}{\langle}
\newcommand{\ran}{\rangle}

\newcommand{\X}{\mathcal{X}}
\DeclareMathOperator{\Irr}{Irr}
\DeclareMathOperator{\adj}{adj}
\DeclareMathOperator{\GL}{GL}
\DeclareMathOperator{\sgn}{sgn}

\newcommand{\dstyle}{\displaystyle}
\DeclareMathOperator{\reg}{reg}
\DeclareMathOperator{\creg}{creg}
\DeclareMathOperator{\sing}{sing}
\DeclareMathOperator{\csing}{csing}
\begin{document}

\title[Character tables]{Submatrices of character tables and basic sets}

\author{Christine Bessenrodt}
\address{Institut f\"{u}r Algebra, Zahlentheorie und Diskrete Mathematik,
Leibniz Universit\"{a}t Hannover,
D-30167 Hannover, Germany}
\email{bessen@math.uni-hannover.de}

\author{J\o rn B.\ Olsson}
\address{
Department of Mathematical Sciences, University of Copenhagen\\
Universi\-tets\-parken~5,
DK-2100 Copenhagen \O, Denmark
}
\email{olsson@math.ku.dk}

\thanks{This work was supported by the Danish Research Council (FNU).}
\subjclass[2010]{Primary 20C30; Secondary 05E10, 05A17}
\keywords{Character tables, basic sets, symmetric groups, determinants,
Cartan matrices, generating functions, $k$-Schur functions.}

\begin{abstract}
In this investigation of character tables of finite groups
we study basic sets and associated representation theoretic
data for complementary sets of conjugacy classes.
For the symmetric groups we find unexpected properties
of characters on restricted sets of conjugacy classes,
like beautiful combinatorial determinant formulae for submatrices
of the character table and Cartan matrices with respect to
basic sets; we observe that similar phenomena occur for the transition
matrices between power sum symmetric functions to bounded partitions
and the $k$-Schur functions defined by Lapointe and Morse. 
Arithmetic properties of the numbers occurring in this context
are studied via generating functions.
\end{abstract}

\date{Version of May 30, 2012}

\maketitle

\section{Introduction}\label{sec:intro}

The investigation of character tables in this paper
draws its motivation and inspiration
from several different sources.

In the $p$-modular representation theory of finite groups,
there has been interest in finding basic sets of characters
on $p$-regular classes. In the more recent investigation
of generalized blocks for the symmetric groups,
also characters on $\ell$-regular classes where $\ell$
is not necessarily  prime have been studied; this is closely
connected to the theory of Hecke algebras at roots of unity. 
For some results on basic sets for finite groups of Lie type see
\cite{Br, G2, G3, GH},
for results on symmetric and alternating groups see
\cite{B-An, BG, BG2, KOR}, \cite[Sec. 6.3]{JK}.

For the symmetric and alternating groups, nice combinatorial
formulae have been found for the determinants
of their character tables and some submatrices of these tables;
via suitable basic sets, this was also connected to the determination
of formulae for the determinant of the Cartan matrices
(see \cite{BO-An, BOS, J, O}).

Another motivation came from the theory of $k$-Schur functions
introduced originally by Lapointe, Lascoux and Morse \cite{LLM};
we are using here the version of the $k$-Schur functions given
by Lapointe and Morse in \cite{LM}. 
The $k$-Schur functions are symmetric functions associated 
to partitions with part
sizes bounded by~$k$, generalizing the classical Schur functions.
An analogue  of the Murnaghan-Nakayama formula was recently
found by Bandlow, Schilling and Zabrocki \cite{BSZ},
which then allowed to compute the transition matrices
between the power sum symmetric functions and the $k$-Schur functions
explicitly.
Indeed, the origin of this paper were observations of the first author
on the determinants of these transition matrices;
this was based on data given in \cite{BSZ}, and further tables
provided by Anne Schilling.
These determinants showed a behavior resembling the one found
for certain submatrices of
character tables of the symmetric groups in \cite{BO-An,BOS,O};
this is described explicitly
in the final section of this paper. 

As the character tables of the symmetric groups are also the transition
matrices between the classical Schur functions and power sum symmetric
functions, with hindsight these analogous phenomena
may not seem to be surprising;
but one should notice that the variant of the
Murnaghan-Nakayama rule presented in \cite{BSZ}
for computing the values of the new transition matrices
for $k$-Schur functions is a rather delicate piece of combinatorics.
Moreover, while the previous investigation of characters of $S_n$
on $p$-regular classes
had been motivated by modular representation theory,
the comparison with $k$-Schur functions now
led us to consider submatrices to bounded partitions
which had not been studied before.
Surprisingly, we found simple combinatorial formulae for the
determinants.
But even more, as in the case of submatrices to
$\ell$-regular classes, the restrictions of characters
to bounded classes give basic sets on these classes,
and we then find simple explicit formulae also for the determinants
of the corresponding Cartan matrices.
For the proof of these properties, it is crucial to
consider also complementary submatrices within the character
table; indeed, we prove here some more general properties
which may be useful also in other circumstances
(especially Theorem~\ref{thm:dualbasic} and its Corollary~\ref{cor:dualbasic}).
In particular, they explain the duality between the
previously observed results for the regular and singular
submatrices in the character table.

\medskip

We now give a brief overview on the following sections.
In Section~\ref{sec:prelim} we consider connections between
complementary submatrices
in square matrices $A$ for which $A^tA$ is diagonal,
with special attention to associated basic sets.
This is then used in Section~\ref{sec:chartables} to investigate
special submatrices of the character tables of
symmetric groups, associated basic sets and corresponding Cartan
matrices. Here, the focus is on determinantal properties,
and the formulae we find involve products of parts or of factorials of
multiplicities of certain partitions.
This motivates to study generating functions for the $p$-powers in such
products in Section~\ref{sec:arithmetic}.
In Section~\ref{sec:applications} we apply our results to
restrictions of characters of the symmetric groups to
regular and singular classes, respectively.
Finally, in Section~\ref{sec:k-Schur} the observations mentioned above
on submatrices of the transition matrices
between power sum symmetric functions to $k$-bounded partitions and
$k$-Schur functions are presented,
and we apply the results of Section~\ref{sec:chartables}
to confirm the validity of some of these observations.

\section{Determinants and basic sets for complementary submatrices}\label{sec:prelim}

An important tool in this section is the Jacobi Minor Theorem which we
restate for the reader's convenience (see \cite{gantmacher} or \cite{P}). 

Let $A=(a_{ij})$ be an $n\times n$ matrix.
We choose  $v$ rows $i_1, \ldots , i_v$ and $v$
columns $k_1, \ldots , k_v$ and let
$M_v$ be the $v$-rowed minor of  $A$
corresponding to this choice of  rows and columns.
By  $M^{(v)}$ we denote the complementary minor to $M_v$,
i.e., the minor to the complementary
rows $i_{v+1}, \ldots, i_n$ and columns $k_{v+1}, \ldots, k_n$
of~$A$.
Corresponding to this numbering we have an associated permutation
$\sigma=\begin{pmatrix} i_1 & \ldots  & i_n \\ k_1 & \ldots & k_n\end{pmatrix}$
which maps $i_j$ to $k_j$, $j=1,\ldots,n$. 

Now let $A'=(A_{ij})$ be the $n \times n$-matrix of cofactors $A_{ij}$
for $A$; this is the transpose of the adjoint matrix to $A$.
Let $M'_v$ be the $v$-rowed minor of $A'$ for the same choice of rows and
columns as for $M_v$. Then the Jacobi Minor Theorem asserts that
$$M'_v =  (\sgn \sigma) \cdot (\det A)^{v-1} M^{(v)} \: .$$

We want to apply this in the following situation which includes
in particular the case of character tables.
\begin{prop}\label{prop:JMT-special}
Let $A\in \GL_n(\C)$ such that
$A^t A =  \Delta$, with $\Delta=\Delta(z_1,\ldots,z_n)$ a diagonal matrix.
Let $A_{(v)}$ be a $v\times v$ submatrix of $A$
for a selection of $v$ rows $i_1, \ldots , i_v$ and $v$
columns $k_1, \ldots , k_v$ in $A$,
$A^{(v)}$ the  submatrix of $A$
to the complementary rows $i_{v+1}, \ldots, i_n$
and columns $k_{v+1}, \ldots, k_n$
of~$A$, and
$\sigma=\begin{pmatrix} i_1 & \ldots  & i_n \\ k_1 & \ldots & k_n\end{pmatrix}$.
Set $\delta_{(v)} =\prod_{j=1}^v z_{k_j}$. 
Then
$$\det A_{(v)} = (\sgn \sigma) \cdot \frac{\delta_{(v)}}{\det A} \, \det A^{(v)} \:.$$
\end{prop}

\proof
We have $(\adj A)^t = (\det A) A \cdot \Delta^{-1}$.
Then the $v$-rowed minor $M'_v$ of this matrix
(to the fixed choice of rows and columns)
is
$$M'_v = (\det A)^v (\det A_{(v)}) (\delta_{(v)})^{-1}\:.$$
By the Jacobi Minor Theorem, we have
$$M'_v = (\sgn \sigma) \cdot (\det A)^{v-1} (\det A^{(v)})\:,
$$
and hence the assertion follows.
\qed

\medskip

The following is an immediate consequence of Cramer's Rule.
\begin{lemma}\label{lem:Cramer}
Let $v, v_1,\ldots,v_k \in \C^k$, $A$ the matrix with rows $v_1,\ldots,v_k$,
$A_i$ the matrix obtained from $A$ by replacing $v_i$ by $v$, $i=1,\ldots,k$.
Suppose that $\det A\neq 0$.
Then $v\in \lan v_1,\ldots,v_k \ran_{\Z }$ if and only if
$\dstyle\frac{\det A_i}{\det A} \in \Z$ for  $i=1,\ldots,k$.
\end{lemma}

Given vectors $v_1,\ldots,v_t$
in $\C^k$, a subfamily
$v_{i_1},\ldots,v_{i_r}$ is defined to be a
{\em basic set} for $v_1,\ldots,v_t$
if it is a $\Z$-basis for
$\lan v_1,\ldots,v_t\ran_{\Z}$.
Lemma~\ref{lem:Cramer} immediately implies a
determinantal criterion for basic sets.
\begin{cor}\label{cor:basicset}
Let $v_1, \ldots, v_t \in \C^k$, $t\geq k$.
Let $A$ be the matrix with rows $v_1,\ldots,v_k$,
$A_{ij}$ the matrix obtained from $A$ by replacing $v_i$ by $v_j$, $i=1,\ldots,k$,
$j=1,\ldots,t$.
Then $v_1, \ldots, v_k$ is a 
basic set for $v_1,\ldots,v_t$
if and only if
$\det A \neq 0$ and
$\dstyle\frac{\det A_{ij}}{\det A} \in \Z$ for all $i,j$. 
\end{cor}

\medskip

When we have a basic set $v_1, \ldots, v_k$ for $v_1, \ldots, v_t$,
the corresponding expansions
$v_i=\sum_{j=1}^k d_{ij} v_j$, $i=1,\ldots,t$,
give an integral {\em decomposition matrix}
$D=(d_{ij})_{1\leq i\leq t \atop 1 \leq j \leq k} \in M_{tk}(\Z)$.
Note that the determinant quotients in the Lemma above
are just these decomposition numbers with
respect to our basic set, by Cramer's Rule.

With this decomposition matrix at hand,
we then call $C=D^tD \in M_k(\Z)$ the corresponding {\em Cartan matrix}.
If we choose a different $\Z$-basis for  $\lan v_1, \ldots, v_t\ran_{\Z }$,
this is related to our basic set by a unimodular transition matrix;
then the Cartan matrix corresponding to this basis is unimodularly
equivalent to~$C$.
As we will only consider invariants such as the determinant or
the Smith normal form, we may thus speak of {\em the} Cartan matrix
for  $\lan v_1, \ldots, v_t\ran_{\Z }$.

Of course, this is closely related to the usual Cartan matrix in
$p$-modular representation theory, which is obtained from the expansion
of ordinary irreducible characters on $p$-regular classes into irreducible
Brauer characters.
It is a classical result that the Brauer characters are a $\Z$-basis
for the ordinary characters on regular classes, but in general
basic sets, i.e.\ subsets of ordinary characters which would give such
a basis on regular classes, are not known.

\medskip

\begin{theorem}\label{thm:dualbasic}
Let $A\in \GL_n(\C)$ with $A^t A = \Delta$, a diagonal matrix.
Let $A_{(v)}$ and $A^{(v)}$ be complementary submatrices of $A$,
corresponding to a selection of $v$ rows and $v$ columns, as in
\ref{prop:JMT-special}, and $\bar A_{(v)}$ and $\bar A^{(v)}$
the $n\times v$ and $n\times (n-v)$
submatrices of $A$ where the rows are restricted to the selected $v$
and $n-v$ column positions, respectively.
Then the rows of $A_{(v)}$ are a basic set for the rows of $\bar A_{(v)}$
if and only if the rows of $A^{(v)}$ are a basic set for the rows of $\bar A^{(v)}$.
\end{theorem}

\proof
We want to use the determinantal criterion for basic sets given above.
For comparing the relevant quotients of determinants,
we apply Proposition~\ref{prop:JMT-special} twice.
We observe that the matrix $\Delta$ does not change when we
interchange rows of~$A$.
Clearly, $\det A_{(v)} \neq 0$ if and only if  $\det A^{(v)} \neq 0$.
Let $A^{ij}_{(v)}$ be the matrix obtained from $A_{(v)}$
by replacing some row $i$ by a row
$j$ in the complementary set,
and $A_{ij}^{(v)}$ the corresponding complementary
submatrix of $A$.
As the factor $\delta_{(v)}$ is the same for both cases, we obtain
$$
\frac{\det A^{ij}_{(v)}}{\det A_{(v)}} = - \frac{\det A_{ij}^{(v)}}{\det A^{(v)}} \:.
$$
Hence the assertion on the basic sets follows by Corollary~\ref{cor:basicset}.
\qed

\begin{cor}\label{cor:dualbasic}
With notation as above,
assume that $(i_1,\ldots,i_n)=(k_1,\ldots,k_n)=(1,2,\ldots,n)$ and
that the rows of $A_{(v)}$ and the rows of $A^{(v)}$ are basic sets for the rows of
$\bar A_{(v)}$ and the rows of $\bar A^{(v)}$, respectively.
Let $d_{ij}, v+1\leq i \leq n, 1 \leq j \leq v$,  be
the corresponding (nontrivial) decomposition numbers
arising from the expansion of the last $n-v$ rows of $\bar A_{(v)}$
w.r.t.\ the rows of $A_{(v)}$,
and let $d'_{ij}, 1\leq i \leq v, v+ 1 \leq j \leq n-v$,
be the (nontrivial) decomposition numbers for~$\bar A^{(v)}$. 
Then these are related by
$$d_{ij} = - d'_{ji} \;, \text{ for } v+1\leq i \leq n, 1\leq j \leq v \:. $$
With $\hat D=(d_{ij})_{v+1\leq i \leq n \atop 1 \leq j \leq v}$,
the Cartan matrices for the two situations are then
$$C_{(v)} = E_v + \hat D^t\hat D \; , C^{(v)} = E_{n-v} + \hat D\hat D^t \:,$$
where $E_m$ is the $m\times m$ identity matrix.
\end{cor}
\proof
The claim on the decomposition numbers follows from the proof of the Theorem.
The assertion on the Cartan matrices is then an immediate consequence.
\qed

\begin{example}\label{S5}
{\rm
For an illustration of the results above,
we consider the character table of the symmetric group~$S_5$.
$$\begin{array}{|r||r|r|r|r|r|r|r|}
\hline
&&&&&&&\\[-5pt]
& (1^5) & (1^3,2) & (1,2^2)
& (1^2,3)    & (2,3) & (1,4)  & (5)\\
&&&&&&&\\[-10pt]
\hline\hline
&&&&&&&\\[-10pt]
[1^5]  & 1 & -1 & 1 & 1 & -1 &  -1 & 1\\
{}[21^3] & 4 & -2 & 0 & 1 & 1 &  0 & -1\\
{}[2^21] & 5 & -1 & 1 & -1 & -1 &  1 & 0\\
{}[31^2] & 6 & 0& -2 &0& 0 &   0 & 1\\
{}[32] & 5 & 1 & 1 & -1& 1 &  -1 & 0\\
{}[41] & 4 & 2 & 0 &1 & -1 &  0 & -1\\
{}[5]     & 1 & 1 & 1 & 1 & 1 & 1 & 1\\ 
\hline
\end{array}
$$
We let $A$ 
be the corresponding matrix, so that $A^tA=\Delta$, the diagonal matrix
with the centralizer orders $120, 12, 8, 6, 6, 4, 5$ on the diagonal.
Let $v=3$, $(i_1,i_2,i_3)=(k_1,k_2,k_3)=(1,2,3)$, so that $A_{(3)}$ is the $3\times 3$
submatrix in the upper left corner and $A^{(3)}$ is the complementary
$4\times 4$ submatrix in the lower right corner of~$A$:
$$A_{(3)}=\begin{pmatrix}
1 & -1 & 1 \\
4 & -2 & 0 \\
5 & -1 & 1
\end{pmatrix}
\quad \text{ and } \quad  A^{(3)}=\begin{pmatrix}
0& 0 &   0 & 1\\
 -1& 1 &  -1 & 0\\
1 & -1 &  0 & -1\\
 1 & 1 & 1 & 1
\end{pmatrix}\:.
$$
The matrix  $\bar A_{(3)}$ is then the submatrix comprising the first
three columns of~$A$,
and $\bar A^{(3)}$ is the complementary submatrix with the last
four columns of~$A$.
Then the theorem says that the rows of  $A_{(3)}$ are a basic set
for the rows of  $\bar A_{(3)}$ if and only if
the corresponding assertion holds for  rows of
$A^{(3)}$ and the rows of  $\bar A^{(3)}$.
Indeed, we will see in Theorem~\ref{thm:small-basicsets}
that this assertion holds for  the chosen
submatrices of the matrix $A$.
In our case, the corresponding decomposition matrices are
(in accordance with the corollary above):
$$
\hat D=\begin{pmatrix}
-3 & 1 &1\\
-1&-1&2\\
-2&-1&2\\
0&-1&1
\end{pmatrix}
\quad \text{, and the dual one } \quad
\begin{pmatrix}
3 & 1 &2&0\\
-1&1&1&1\\
-1&-2&-2&-1
\end{pmatrix}\:.
$$
The corresponding Cartan matrices are
$$C_{(3)}=\begin{pmatrix}
15 & 0 & -9 \\
0 & 5 & -4 \\
-9 & -4 & 11
\end{pmatrix}
\quad \text{ and } \quad
C^{(3)}=\begin{pmatrix}
12& 4 &   7 & 0\\
 4& 7 &  7 & 3\\
7 & 7 &  10 & 3\\
 0&3&3&3
\end{pmatrix}\:.
$$
We will come back to these matrices in section~\ref{sec:chartables}
when we study submatrices of the character tables of symmetric groups
and certain Cartan matrices associated to them in detail.
}
\end{example}

\section{Submatrices of the character table of $S_n$}\label{sec:chartables}

With the preparations in the previous section
done for a more general situation, we now want to apply
this in the context of character tables of symmetric groups.
For the background on the representation theory of the symmetric groups,
the reader is referred to \cite{JK};
for the connection to symmetric functions see
also  \cite{M,Sagan,EC2}.

\smallskip

In the theory of symmetric functions, the restriction to $k$-bounded partitions
(i.e., those with largest part at most $k$) has led to the notion
of $k$-Schur functions; letting $k$ increase gives a filtration
of the algebra of symmetric functions with nice properties
(see \cite{BSZ, LM}).
As mentioned in the introduction, observations on the
transition matrices between these
and the power sum functions led to the investigations
on character tables of symmetric groups pursued here.
While this started by considering submatrices to $k$-bounded
partitions, it then appeared that indeed
more refined properties hold for more general subtables.
This is what we want to present in this section.

\medskip

We fix a positive integer $n$,
and we briefly recall some notation.

Let $P(n)$ be the set of all partitions of $n$.
If $\mu\in P(n)$,
then $z_{\mu}$ denotes the order of the centralizer of an element of
cycle type $\mu$ in $S_n$.
Explicitly, when $\mu$ is written in exponential notation,
i.e., $\mu=(1^{m_1(\mu)}, 2^{m_2(\mu)},\ldots)$,
we have $z_\mu=a_{\mu}b_{\mu}$,  where
$$a_{\mu}=\prod_{i \ge 1}i^{m_i(\mu)}, \;
b_{\mu}=\prod_{i \ge 1}m_i(\mu)!\:.$$

We will always order the partitions of $n$ lexicographically,
where this means lexicographic reverse order when the parts of the
partitions are written in increasing order.
We denote this relation
by $\la\leq \mu$ or $\la < \mu$, if $\la\neq \mu$.
For example, for $n=4$,
the ordered list is
$(1^4) < (1^2,2) < (2^2) < (1,3) < (4)$.

For $\la \in P(n)$,
we denote by $\chi^{\la }$ the corresponding
irreducible character of $S_n$.
For $\mu\in P(n)$, $\chi^{\la }_\mu = \chi^{\la }(\sigma_\mu)$ then denotes the
character value of $\chi^{\la }$ on
an element $\sigma_\mu \in S_n$ of cycle type $\mu$.

We let $X=(\chi^{\la }_\mu )_{\la,\mu\in P(n)}$ be the character table
of $S_n$. This is also the transition matrix between the
classical Schur functions and the power sum symmetric functions.

\medskip

For a given $\al\in P(n)$, let
$X_{(\al)} = (\chi^{\la }_\mu )_{\la,\mu <\al}$
and $X^{(\al)} = (\chi^{\la }_\mu )_{\la,\mu \geq \al}$
be the submatrices corresponding to the ``small'' and ``large''
partitions with respect to $\al$, respectively.
Note that the only reason for defining $X_{(\al)}$ this way, rather
than running over all $\la, \mu \leq \al$, 
is that our focus will be on
complementary submatrices; in particular, adjoining a further maximal
element $\hat 1$ to $P(n)$ will also give the case where $X_{(\hat 1)}=X$.

For the submatrix  $X^{(\al)}$ whose entries are the values of characters
labelled by  ``large'' partitions on classes of ``large'' cycle type 
we may now easily determine
the corresponding Smith normal form.
For an integer matrix $A$ we denote  its Smith normal form by $\S(A)$,
and for  integers $x_1, \ldots , x_n$ we denote by
$\S(x_1, \ldots, x_n)$ the Smith normal form of the
diagonal matrix with entries $x_1, \ldots , x_n$ on the diagonal.
\smallskip

In the proof of the following result we will use the connection
of the irreducible characters of the symmetric groups
to the permutation characters.
Here,  the permutation character $\xi^{\la}$ associated to 
the partition $\la$ is obtained by inducing the trivial character of the Young
subgroup $S_{\la}$ up to $S_n$.
The values of the permutation characters appear
in the expansion of the power sum symmetric functions
into monomials.

\begin{theorem}\label{thm:smith}
Let $\al\in P(n)$. Then
$$\det X^{(\al)} = \prod_{\mu \geq  \al} b_\mu $$
and
$$\S(X^{(\al)}) = \S(b_\mu; \mu \geq  \al) \:.$$
\end{theorem}

\proof 
The transition matrix between the permutation characters and the
irreducible characters of $S_n$ is an upper unitriangular matrix (with respect to
our chosen order), see \cite[Sec.~2.2]{JK}.
Thus $X^{(\al)}$ is unimodularly equivalent to the corresponding
permutation character matrix
$\Xi^{(\al)} = (\xi^{\la}(\mu))_{\la,\mu \geq  \al}$.
By \cite[Corollary 11]{BOS}, this is a lower triangular matrix
with entries $b_{\mu}$, $\mu \geq \al$, on the diagonal,
where each $b_\mu$ divides all entries in the same row.
Thus the statement on the determinant follows; also,
$\S(\Xi^{(\al)})= \S(b_\mu;\mu \in P(n), \mu \geq  \al)$,
yielding the claim on the Smith normal form for~$X^{(\al)}$.
\qed

\medskip

The reasoning above also has the following consequence. 
\begin{cor}\label{cor:large-basicsets}
Let $\al \in P(n)$.
Then both the
permutation characters $\xi^{\la}$, $\la \in P(n)$, $\la \geq  \al$,
as well as the irreducible characters $\chi^{\la}$, $\la \in P(n)$, $\la \geq  \al$,
provide basic sets for the characters restricted to classes of
cycle type $\geq \al$.
\end{cor} 
\proof
As the permutation character matrix $\Xi$ is a lower triangular matrix,
the rows of any $k\times k$ principal submatrix in the lower right corner
give a basic set for the rows of the submatrix of $\Xi$ comprising the final
$k$ columns.
Since the transition matrix between the permutation characters and the
irreducible characters of $S_n$ is upper unitriangular, this proves both that
the rows of $\Xi^{(\al)}$ as well as that the rows of $X^{(\al)}$ are
a basic set for the characters on classes of
cycle type $\geq \al$.
\qed
\medskip

As the character table satisfies
$$X^tX=\Delta(z_\mu; \mu \in P(n))$$
(with respect to our chosen order),
we can now apply the results on complementary submatrices
from the previous section to obtain corresponding results also for
the characters labeled by ``small'' partitions
on classes of ``small'' type.

In fact, the following result was suspected after having observed
on the data in \cite{BSZ} similar behavior of the transition matrices
for the $k$-Schur functions.

\begin{theorem}\label{thm:small-basicsets}
Let $\al \in P(n)$.
Then the irreducible characters $\chi^{\la}$, $\la \in P(n)$, $\la < \al$,
provide a basic set for the characters restricted to classes of cycle type
$< \al$.

Furthermore,
$$\det X_{(\al)} = \prod_{\mu < \al} a_\mu \:.$$
\end{theorem}

\proof
The first statement follows immediately from Corollary~\ref{cor:large-basicsets}
using Theorem~\ref{thm:dualbasic}.

Towards the second assertion, Proposition~\ref{prop:JMT-special} gives
\begin{align*}
\det X_{(\al)} &=
\frac{\prod_{\mu <  \al} z_\mu}{\det X} \det X^{(\al)} \\
&=
\frac{\prod_{\mu <  \al} z_\mu \prod_{\mu \geq  \al} b_\mu}{\det X} \:.
\end{align*} 
Recall that $z_\mu=a_\mu b_\mu$;
now,  by \cite[Cor. 6.5]{J} and Theorem~\ref{thm:smith}
(see also \cite{O} and~\cite{SS}) we have
$$\det X = \prod_{\mu \in P(n)} a_\mu = \prod_{\mu \in P(n)} b_\mu  $$
and hence the claim follows.
\qed

\medskip

\begin{rems}{\rm
(1)
The two determinant formulae for $X_{(\al)}$ and $X^{(\al)}$ together
may be viewed as giving a very nice
interpolation between the two expressions for
the determinant of the full character table that we have also used above, i.e.,
$\det X = \prod_{\mu \in P(n)} b_\mu$ and $\det X = \prod_{\mu \in P(n)} a_\mu $.

(2)
One might suspect for the Smith normal form that
$\S(X_{(\al)}) = \S(a_\mu,\mu <  \al)$,
but this is not true as one already sees in small character tables,
for example, taking the full character table of $S_4$. 
}
\end{rems}

\begin{example}
{\rm
In Example~\ref{S5}, we have considered for the partition $\al=(1^2,3)$
the submatrices $X_{(\al)}$ and
$X^{(\al)}$ of the character table of $S_5$,
which appeared there as $A_{(3)}$ and $A^{(3)}$, respectively.
We had already pointed out the basic set property in Example~\ref{S5}.
Illustrating the results above, we have
$$\det(X_{(1^2,3)})= 8= a_{(1^5)} \cdot a_{(1^3,2)} \cdot  a_{(1,2^2)} 
\;,$$
$$\det(X^{(1^2,3)})=2=
b_{(1^2,3)}\cdot b_{(2,3)}\cdot b_{(1,4)}\cdot b_{(5)} 
,\:$$
$$
\S(X^{(1^2,3)})=\Delta(1,1,1,2)=\S(b_{(1^2,3)}, b_{(2,3)},b_{(1,4)},b_{(5)})\:.$$
Here, $\S(X_{(1^2,3)})=\Delta(1,2,4)=
\S(a_{(1^5)}, a_{(1^3,2)}, a_{1,2^2)})$ does hold,
but as remarked above, this is  not true in general.
}
\end{example}

Next we want to take a look at the Cartan matrices corresponding to
the basic sets given above.
First, we state a more general easy Lemma.
\begin{lemma}\label{lem:Cartandet}
Let $G$ be a finite group, $\Irr(G)=\{\chi_1,\ldots,\chi_t\}$.
Select conjugacy classes $C_1,\ldots,C_k$, and denote by $\chi'$
the restriction of the character $\chi$ to the union of
these classes; let $z_i=|G|/|C_i|$ for $i=1,\ldots,k$.
Assume that $\{\chi_1', \ldots,\chi_k'\}$ is a basis for
$\lan\chi_1',\ldots,\chi_t'\ran_{\C }$.
Let $Y=(\chi_i(C_j))_{1\leq i,j\leq k}$ and
$\bar X=(\chi_i(C_j))_{1\leq i\leq t\atop 1\leq j\leq k}$ be
corresponding submatrices of the character table,
$D=\bar X Y^{-1}$ the associated decomposition matrix
(not necessarily integral),
and $C=D^tD$.

Then $$\det C = \dstyle \frac{\prod_{i=1}^k z_i}{(\det Y)^2}\:.$$
\end{lemma}

\proof
By the orthogonality relations for characters we have
$$C=D^tD=(Y^{-1})^t {\bar X}^t \bar X Y^{-1} =
(Y^{-1})^t \Delta(z_i; i=1,\ldots,k) Y^{-1}\:,$$
and hence the claim follows.
\qed

\medskip

We now return  to the context of the character table
of the symmetric group $S_n$ studied before
and consider the associated Cartan matrices.
\begin{theorem}\label{thm:Cartan-det}
Let $\al \in P(n)$, and let $C_{(\al)}$ (or $C^{(\al)}$, respectively)
be the Cartan matrix corresponding to the basic set of character restrictions
associated to the partitions $<  \al$ (or $\geq  \al$, respectively).
Let $a_{(\al)}= \prod_{\mu <  \al} a_\mu$,
$b_{(\al)}= \prod_{\mu <  \al} b_\mu$,
and let  $a^{(\al)}, b^{(\al)}$ be the complementary products.
Then we have
$$\det C_{(\al)} = \frac{b_{(\al)}}{a_{(\al)}}
 = \frac{a^{(\al)}}{b^{(\al)}} = \det C^{(\al)}\:.$$
\end{theorem}

\proof
Let $z_{(\al)}= \prod_{\mu <  \al} z_\mu$.
Using Lemma~\ref{lem:Cartandet} and Theorem~\ref{thm:small-basicsets},
we have
$$\det C_{(\al)} = \frac{z_{(\al)}}{a_{(\al)}^2}=\frac{b_{(\al)}}{a_{(\al)}}\:.
$$
Similarly, using Theorem~\ref{thm:smith}, we obtain the
formula for $\det C^{(\al)}$.
But as
$$a_{(\al)} \cdot a^{(\al)} =
\prod_{\mu \in P(n)} a_\mu = \prod_{\mu \in P(n)} b_\mu
= b_{(\al)} \cdot b^{(\al)} \:,$$
we also get the equality in the middle.
\qed

\medskip

\begin{example}
{\rm
We take another look at Example~\ref{S5}.
As said above,
the selection of rows and columns chosen there corresponds
in the notation introduced here to the case $\al=(1^2,3)$,
i.e., the Cartan matrices $C_{(3)}$ computed there are
the Cartan matrices $C_{(1^2,3)}$ and $C^{(1^2,3)}$ here.
Indeed, we find for the determinants
$$\det C_{(1^2,3)} = \frac{5! \cdot 3! \cdot 2!}{1\cdot 2 \cdot 4}
= 180 = \frac{3\cdot 6 \cdot 4 \cdot 5}{2\cdot 1 \cdot 1 \cdot 1} = \det C^{(1^2,3)}\:.$$
}
\end{example}

\medskip

As a consequence of the above, we have in particular the following
interesting arithmetic result.
\begin{cor}\label{cor:integralquotients}
Let $\al \in P(n)$.
Then the quotient
$\dstyle \frac{b_{(\al)}}{a_{(\al)}}=\frac{a^{(\al)}}{b^{(\al)}}$
is an integer.
\end{cor}

\proof
For $\al \in P(n)$, the Cartan matrix $C_{(\al)}$
is an integral matrix, as it is associated to a
basic set of characters.
Hence its determinant is an integer,
and thus by Theorem~\ref{thm:Cartan-det}
the quotient
$\frac{b_{(\al)}}{a_{(\al)}}=\frac{a^{(\al)}}{b^{(\al)}}$
is an integer.
\qed

\medskip

It is natural to ask whether there is also a more direct
combinatorial proof of the arithmetic property above.
Indeed, in Section~\ref{sec:arithmetic} we will discuss in more detail
the case $\al=(1^{n-k-1},k+1)$, i.e., the arithmetics in the situations
where we take in the quotient $\frac{b_{(\al)}}{a_{(\al)}}$
products over all $k$-bounded partitions. We have not been able to
generalize the arguments to the case of general~$\al$.

\section{The arithmetic of the $a$- and $b$-numbers}\label{sec:arithmetic}

In this section we study
generating functions for the
$p$-part of products of $a_\mu$'s for certain sets of partitions $\mu$ ($a$-numbers),
as well as in the products of  $b_\mu$'s for certain sets
of partitions $\mu$ ($b$-numbers), where $p$ is a prime number.

For any partition $\mu$  we had defined the numbers
$a_{\mu}$ and $b_{\mu}$ in Section~\ref{sec:chartables}.
For any set of partitions $Q$
we define
\beq
a_{Q}=\prod_{\mu \in Q}a_{\mu}, \hspace{.1 in}
b_{Q}=\prod_{\mu \in Q}b_{\mu}.
\eeq
Note that we have already used in Section~\ref{sec:chartables}
the fact that
$$
a_{P(n)}=b_{P(n)} \quad \text{for all } n \in \N \:.
$$
For the set $P(n,k)$ of partitions of $n$ with
largest part at most $k$, the corresponding equality does not hold
in general, as can already been seen in the case $k=1$, $n>1$.

Using generating functions we compute
for a fixed $k$ and a given
prime number~$p$
the exponents of $p$ dividing
$a_{P(n,k)}$ and $b_{P(n,k)}$.
In particular, we show that they
are the same in the expressions
$a_{P(n)}$ and $b_{P(n)}$;
this gives yet another proof of the equality above. 

For a given set $Q$ of partitions the numbers $a_Q$ and $b_Q$ are defined as
products of some lists of numbers. Rather than studying the valuation of each
prime number in $a_Q$ and $b_Q$, one could look at the multiplicity of each
of the integers in these lists and try to prove results about these
multiplicities. However examples show that it is not possible to prove
a result like Proposition~\ref{prop:abS} below in this way.

\smallskip

Let $p$ be a  fixed prime; for $n \in \N$ we denote by
$\nu_p(n)$ the exponent of the maximal $p$-power
dividing $n$, and we then define the $p'$-part $w(n)$ of $n$ by the
equation 
\beq\label{uv}
n=p^{\nu_p(n)}w(n).
\eeq

Furthermore, it will be convenient to treat some situations in greater generality
and consider partitions with parts from a selected set $S\subseteq \N$.
Then $P(n,S)$ denotes the number of partitions of $n$ with all parts 
being in~$S$.
For example, $S=\N$, $S=\{1,\ldots,k\}$ or $S=\{j\in \N \mid \ell \nmid j\}$
are natural choices, giving all partitions, $k$-bounded partitions or $\ell$-regular partitions,
respectively. When $S$ is not mentioned explicitly, we understand this to be $S=\N$.

\medskip

We note the following easy facts on generating functions involving $S$.
The generating function for the number of partitions in $P(n,S)$ is
$$P_S(q)=\prod_{j\in S} \frac{1}{1-q^j}\:.$$
Now we fix $i\in S$ and $m\in \N$.
Then the generating function for the number of partitions
in $P(n,S)$  with at least $m$ parts $i$ is
$$P_S(q)q^{im}\:, $$
and hence the generating function for the number of partitions
in $P(n,S)$  with exactly $m$ parts $i$ is
$$P_S(q)(1-q^i)q^{im}\:.$$
Thus the generating function for the total number of parts $i$ in all partitions
in $P(n,S)$  is
$$P_S(q)(1-q^i)\sum_{m\ge 1}mq^{im}= P_S(q)\frac{q^i}{1-q^i}\:.$$

\medskip

As we will see, generating functions for the number of divisors
of $n \in \N$ will play a special role.
Let $t_S(n)$ be the number of divisors $d\in S$ of $n$.
The generating function $T_S(q)$ for $t_S(n)$ is then
\beq
  T_S(q)=\sum_{i \in S}\frac{q^i}{1-q^i} \:.
\eeq
For $r\in \N$, then
$T_S(q^r)=T_{rS}(q)$ is the generating function for the number of
divisors $rd$ of $n$ with $d\in S$.
In particular, $T(q^{p^v})=T_{\N }(q^{p^v})$ is the generating
function for the number of divisors $d$ of $n$ with $\nu_p(d)\ge v$.
For $S=\{j\in \N \mid p\nmid j\}$, we also write
$T_p(q)$ for the generating function $T_S(q)$.
Thus $T_p(q^{p^v})$ is the generating function for the number
of divisors $d$ of $n$ with $\nu_p(d)=v$.

We have the following connection between $T_S(q)$ and $P_S(q)$
which is immediate from the observations made on these functions above;
note that special cases (with basically the same proof) are contained
in~\cite{BO-Cartan}.
Here, we denote by $\ell(\la)$ the
length of the partition~$\la$,
i.e., the number of (positive) parts of~$\la$. 
\begin{prop}\label{prop:LPT}
Let $\ell_S(n)=\sum_{\la \in P(n,S)} \ell(\la)$ and
$L_S(q)$ the corresponding generating function.
Then
$$L_S(q)=P_S(q)T_S(q) \:.$$
\end{prop}

\medskip

Next we want to consider generating functions for
weights on the parts and divisors, respectively,
according to their $p$-value.
The generating function for the total $p$-weight
$$e_{S,p}(n)=\sum_{d\mid n \atop d \in S} \nu_p(d)$$
of the divisors of $n$ in our given set $S$ is
\beq
E_{S,p}(q)= \sum_{n\geq 1} e_{S,p}(n)q^n = \sum_{r \in S} \nu_p(r)\frac{q^{r}}{1-q^{r}} \:.
\eeq
We consider also the closely related generating function
\beq
F_{S,p}(q) = \sum_{n\geq 1}f_{S,p}(n)q^n
:= \sum_{r\in S}\sum_{j \ge 1}\frac{q^{rp^j}}{1-q^{rp^j}} \:.
\eeq

\medskip

For the situations we are interested in here, two properties of
the subset $S$ of $\N$ are particularly important.
We call $S$ {\em $p$-divisible} if $pr \in S$, for some $r\ge 1$,
implies $r\in S$.
The set $S$ is called {\em $p$-closed} if for any  $d\in S$
also $pd \in S$.

\begin{rem}
{\rm
Note that for $S=\N$,
each number $d\in S$ can be written
in $\nu_p(d)$ ways in the form $p^jr$, with $r\in S$,
$j\ge 1$. Hence in this case we have an equality between
the two generating functions just defined, i.e.,
$E_{S,p}(q)=F_{S,p}(q)$ for~$S=\N$.
}
\end{rem}

\begin{prop}\label{prop:F-coeff}
Let $S\subseteq \N$ be $p$-divisible.
Then for all $n\in \N$ we have
$$f_{S,p}(n)=\nu_p(n) t_S(n) - e_{S,p}(n) $$
and
$$f_{S,p}(n) \geq e_{S,p}(n) \:.$$

If $S$ is in addition $p$-closed (in particular,
when $S=\N$), we have
$$e_{S,p}(n)=f_{S,p}(n)=  {\nu_p(n)+1 \choose 2} \ t_S(w(n)) \:,$$ 
and hence 
$$E_{S,p}(q)=F_{S,p}(q)=
\sum_{v \ge 1} \binom{v+1}{2}T_S(q^{p^{v}}) \:.$$ 
\end{prop}

\proof
By definition, $f_{S,p}(n)=|\{(r,j)\mid r\in S, j\geq 1, rp^j\mid n\}|$.
We want to consider the contribution coming from a maximal
$p$-string $u, pu, \ldots, p^cu$ of divisors of $n$ in~$S$, where
$u\mid w(n)$, and clearly $c\leq \nu_p(n)$ (for the form of the maximal string
we have used the assumption that $S$ is $p$-divisible).
Any divisor $r=p^iu$ in this string contributes the
pairs $(r,p^j)$, $j=1,\ldots,\nu_p(n)-i$, to the count;
hence from the complete string we have the contribution
$$\sum_{i=0}^c (\nu_p(n)-i) = \nu_p(n) (c+1) - \frac{c(c+1)}2\:.$$
Note also that the contribution of the string to the coefficient
$e_{S,p}(n)$ is $\sum_{i=0}^c i$, since $ \nu_p(p^iu)=i$; 
as $c\le \nu_p(n)$, this is
less than or equal to the contribution of the string for $f_{S,p}(n)$,
and equality holds only when $c=\nu_p(n)$.
The latter condition is always satisfied when $S$ is $p$-closed.

Thus, we get in total
$$f_{S,p}(n)=\sum_{r\in S \atop r \mid n} (\nu_p(n)-\nu_p(r))=\nu_p(n) t_S(n)-e_{S,p}(n)\:,$$
and we have also shown that the claimed inequality holds.

Furthermore, when $S$ is in addition $p$-closed,
all $p$-strings start with a divisor of $w(n)$ and then
are of full length $c+1=\nu_p(n)+1$, and thus each one gives a contribution
${\nu_p(n)+1 \choose 2}$, and the total count for $f_{S,p}(n)$ is as claimed.
\qed

\bigskip

After these preparations, we can now move on to compute the
generating functions $A_{S,p}(q)$ and $B_{S,p}(q)$ for the
exponents $a_{S,p}(n)$ and $b_{S,p}(n)$ of $p$
dividing $a_{P(n,S)}$ and $b_{P(n,S)}$, respectively.

\begin{prop}\label{ab-gf}
Let $S\subseteq \N$. Then we have
$$ A_{S,p}(q)=P_S(q)E_{S,p}(q) \: , \:
B_{S,p}(q)=P_S(q)F_{S,p}(q) \:.
$$
\end{prop}

\proof
By definition and the formula for the generating function for the total
number of parts $i$ in all partitions in $P(n,S)$, we obtain
$$
A_{S,p}(q)= \sum_{r\in S} \nu_p(r) P_S(q) \frac{q^r}{1-q^r}
= P_S(q) E_{S,p}(q) \:.
$$

We next consider the generating function $B_{S,p}^{(i)}(q)$ for the exponent
of $p$ in the factorials $m!$ coming from parts $i\in S$
in partitions of $P(n,S)$ with multiplicity $m$.
The exponent of $p$ in $m!$ is
\beq
\nu_p(m!)=
\left\lfloor \frac{m}{p}\right\rfloor +\left\lfloor \frac{m}{p^2}\right\rfloor + \ldots
= |\{(s,t)| s, t \ge 1,  p^st \le m\}| \:.
\eeq
We thus get a contribution $1$ to the exponent whenever
$m\ge p^st$, for some $s,t \ge 1$.
Hence we obtain
\beq
B_{S,p}^{(i)}(q)=\sum_{s,t \ge 1} P_S(q)q^{ip^st}=
P_S(q)\sum_{s \ge 1}\frac{q^{p^si}}{1-q^{p^si}}\:.
\eeq
Summing over all $i \in S$ we get
$$
B_{S,p}(q)=\sum_{i\in S} B_{S,p}^{(i)}(q)=
P_S(q)\sum_{i\in S, s \ge 1}\frac{q^{p^si}}{1-q^{p^si}}=
P_S(q)F_{S,p}(q)\:.
\qedhere
$$

\medskip

\begin{prop} \label{prop:abS}
Let $S\subseteq \N$ be $p$-divisible. Then for all $n\in \N$ we have
$$\nu_p(a_{P(n,S)}) \le  \nu_p(b_{P(n,S)}) \:,$$
and equality holds if $S$ is also $p$-closed.
\end{prop}

\proof
Let $p$ be a prime;
by Proposition~\ref{ab-gf} and Proposition~\ref{prop:F-coeff} we have
$$
a_{S,p}(n)=\sum_{r=0}^n p_S(n-r)e_{S,p}(r) \le
\sum_{r=0}^n p_S(n-r)f_{S,p}(r)=b_{S,p}(n)\:,
$$
and equality holds if $S$ is also $p$-closed.
Hence the exponent of $p$ in  $a_{P(n,S)}$ is less than or equal
to  the exponent of $p$ in  $b_{P(n,S)}$.
\qed

\medskip

We immediately deduce the following result
which provides in particular
the alternative proof for $a_{P(n)}=b_{P(n)}$ mentioned earlier.
\begin{cor} \label{cor:abS}
Let $S\subseteq \N$ be $p$-divisible for all primes $p$.
Then for all $n\in \N$ we have
$$a_{P(n,S)} \mid  b_{P(n,S)} \:.$$
If $S$ is in addition $p$-closed for all primes~$p$, we have
equality for all $n\in \N$, i.e.,
$$A_{S,p}(q)=B_{S,p}(q)\:.$$
In particular,
$$A_p(q)=B_p(q)\:.$$
\end{cor}

As $S=\{1,\ldots,k\}$ is $p$-divisible for all primes $p$,
we may also deduce the divisibility property obtained in
Corollary~\ref{cor:integralquotients} in our more special
situation by very different means.
\begin{cor} \label{abk}
Let $k\in \N$. Then for all $n\in \N$ we have
$$a_{P(n,k)} \mid b_{P(n,k)} \:.$$
\end{cor}

\section{Applications to regular and singular character tables}\label{sec:applications}

The results on complementary submatrices
give further insights for the regular character tables
and their singular counterparts,
studied in earlier papers.

We fix a natural number $\ell$, which need not be a prime. 
In \cite{KOR} and \cite{O},
the characters on $\ell$-regular conjugacy classes of $S_n$ are investigated.
We recall the notation used in these papers.
A conjugacy class of cycle type $\mu$ is $\ell$-regular if no part of $\mu$
is divisible by~$\ell$, and we then call $\mu$ an $\ell$-class regular partition;
otherwise the class is called $\ell$-singular, and the partition
is then $\ell$-class singular. 
The restriction of a character $\chi$ to the $\ell$-regular classes
is denoted by $\chi_{(\reg)}$. 
Recall that a partition is said to be $\ell$-regular if no part size appears
with multiplicity $\geq \ell$; otherwise it
is called $\ell$-singular. 
The $\ell$-regular character table $X_n^{\reg}$ is then the submatrix to the
characters
$\chi^{\la}$, $\la$ $\ell$-regular, on the $\ell$-regular classes.
Finally, we set \,
$\dstyle a_n^{\creg}= \prod_{\mu \in P(n) \atop \ell\text{-class regular}} a_\mu$.
Similarly for the product of $b_\mu$'s and $\ell$-class singular partitions.

We recall the following results from \cite{BOS, KOR, O}. 
\begin{prop}
Let $n\in \N$.
\begin{enumerate}
\item[{(1)}] \cite[Prop. 4.2]{KOR}
The character restrictions $\chi_{(\reg)}^\la$, $\la\in P(n)$ $\ell$-regular,
form a basic set for the character restrictions on $\ell$-regular
classes. 
\item[{(2)}] \cite[Theorem 2]{O}
For the regular character table, we have $|\det X_n^{\reg}| = a_n^{\creg}$.
\item[{(3)}] \cite{BOS, O}
For the Cartan matrix $C^{\reg}$ to the $\ell$-regular classes and characters,
we have \, $\dstyle \det C^{\reg} = \frac{b_n^{\creg}}{a_n^{\creg}}$.
\end{enumerate}
\end{prop}

Using our results on complementary submatrices, we immediately have the following
results on characters associated
to $\ell$-singular partitions restricted to $\ell$-singular classes;
in fact, for part (2) below a similar proof was given in~\cite{BOS}. 

\begin{cor}
Let $n\in \N$.
\begin{enumerate}
\item[{(1)}]
The character restrictions $\chi_{(\sing)}\la$, $\la\in P(n)$ $\ell$-singular,
form a basic set for the character restrictions on $\ell$-singular
classes.
\item[{(2)}] \cite[Theorem 3]{O}
For the singular character table, we have $|\det X_n^{\sing}| = b_n^{\csing}$.
\item[{(3)}]
For the Cartan matrix $C^{\sing}$ to the $\ell$-singular classes and characters,
we have
$$\dstyle \det C^{\sing} = \frac{a_n^{\csing}}{b_n^{\csing}}
= \frac{b_n^{\creg}}{a_n^{\creg}}
= \det C^{\reg}\:.$$
\end{enumerate}
\end{cor}

\medskip

The set of $\ell$-regular partitions of $n$ is just $P(n,S)$ with
$S=\{m\in \N \mid \ell \nmid m\}$.
This set $S$ is clearly $p$-divisible for all primes~$p$.
Thus we can also apply some of the results on the arithmetic of the $a$- and $b$-numbers.
In particular, Proposition~\ref{prop:abS} also gives a combinatorial explanation for the property
$$a_n^{\creg}\mid b_n^{\creg}$$
that is implied above by the fact that the quotient is the determinant
of an integral matrix.
Furthermore, as the set $S$ is also $p$-closed for all primes $p$ not dividing $\ell$,
 by Proposition~\ref{prop:abS} the quotient
$\dstyle\frac{b_n^{\creg}}{a_n^{\creg}} \in \Z$
is only divisible by primes dividing $\ell$.
In fact, it has been shown in \cite{O} that this number is a power of~$\ell$.

\section{On $k$-Schur functions: Observations and applications}\label{sec:k-Schur}

As mentioned in the introduction, the initial motivation for this
investigation of restricted character tables came from
observations on the expansion coefficients of power sum
symmetric functions to $k$-bounded partitions into
$k$-Schur functions.
For more precise statements, we briefly recall
some notation from \cite{BSZ}.

We fix an integer $n\in \N$ and consider only partitions of $n$ in this section. 
We let $ P^{(k)}$
denote the set of $k$-bounded partitions;
for $\la\in  P^{(k)}$,   $s_{\la}^{(k)}$ denotes the corresponding
$k$-Schur function as defined by Lapointe and Morse in~\cite{LM}. 
The set of $k$-Schur functions $s_{\la}^{(k)}$, $\la \in P^{(k)}$,
forms a basis for the space spanned by the homogeneous symmetric functions
$h_{\la}$, $\la\in P^{(k)}$ \cite[Property 27]{LM}.
Hence, for $\nu\in P^{(k)}$ the corresponding power sum symmetric function
can be decomposed as
$$p_\nu =
\sum_{\la \in P^{(k)}} \chi_{\la,\nu}^{(k)} \, s_{\la}^{(k)}\:.
$$
By the work of Bandlow, Schilling and Zabrocki \cite{BSZ},
the coefficients $\chi_{\la,\nu}^{(k)}$ can be computed combinatorially by
an intricate analogue of the Murnaghan-Nakayama formula involving
$k$-ribbon tableaux (see \cite{BSZ} for details).
As the $k$-Schur functions are not self-dual, one may also  consider the
coefficients $\tilde\chi_{\nu,\la}^{(k)}$
in the expansion of the power sums into the dual $k$-Schur functions;
these also appear in the expansion
$$s_\nu^{(k)} =
\sum_{\la \in P^{(k)}} \frac{1}{z_\la}\tilde\chi_{\nu,\la}^{(k)}\, p_{\la}\:.
$$
With $\X^{(k)}=(\chi_{\la,\nu}^{(k)})_{\la,\nu \in P^{(k)}}$ and
$\tilde \X^{(k)}= (\tilde\chi_{\la,\nu}^{(k)})_{\la,\nu \in P^{(k)}}$
(with the usual order),
the duality implies that we have for the product
$$(\X^{(k)})^t\cdot \tilde \X^{(k)}
=\Delta(z_\la, \la \in P^{(k)})\:.$$
Thus the table $\tilde \X^{(k)}$ can be computed from the table $\X^{(k)}$,
but a  combinatorial formula for the coefficients $\tilde\chi_{\la,\nu}^{(k)}$
has not yet been obtained.

In the available data
(tables in \cite{BSZ} and further tables provided by Anne Schilling),
we have discovered close connections between $\X^{(k)}$,
$\tilde\X^{(k)}$ and the
restricted character table
$X^{(k)}:=(\chi^{\la }_\mu )_{\la,\mu \in P^{(k)}}$ 
corresponding to the $k$-bounded partitions
(beware not to confuse this with the matrix $X^{((k))}$).

As we have done it for the character table,
we may also take a closer look at submatrices of the tables
$\X^{(k)}$ and $\tilde\X^{(k)}$.
For a given $\al \in P^{(k)}$,
the corresponding (upper) principal submatrix
$$\X^{(k)}_{(\al)}=(\chi_{\la,\nu}^{(k)})_{\la,\nu < \al}$$
of $\X^{(k)}$ is closely related to the corresponding (upper) principal
submatrix $X_{(\al)}$ of the character table,
refining the observations mentioned above.

Indeed, based on the properties of the $k$-Schur functions
studied in \cite{LM} and our results in Section~\ref{sec:chartables},
we can prove the following.

\begin{theorem}\label{thm:k-Schur}
Let $k\in \N$, $\al \in P^{(k)}$.
\begin{enumerate}
\item
The matrices $X^{(k)}$ and $\X^{(k)}$, and
the matrices $X_{(\al)}$ and $\X^{(k)}_{(\al)}$, respectively,
are related by integral lower unitriangular
transition matrices.

\item For the determinants we have
$$
\det \X^{(k)} = \prod_{\la \in P^{(k)}} a_{\la} = \det X^{(k)}
\; , \;
\det(\X^{(k)}_{(\al)})=
\prod_{\la < \al} a_{\la} =\det(X_{(\al)})\:.
$$
\end{enumerate}
\end{theorem}
\proof
(1) Let $\rhd$ denote the dominance order on partitions. 
As observed in \cite[Property 28]{LM}, for any $\la\in P^{(k)}$ we have
$$
s_{\la}^{(k)}= s_\la + \sum_{\mu :\ \mu \rhd \la} d_{\la\mu}^{(k)} s_\mu\;, \text{ for }  d_{\la\mu}^{(k)}\in \Z.
$$
Hence for $\la, \nu \in P^{(k)}$ we have
by the expansion formulae for the power sums
$$
\chi^{\la}_{\nu}=\sum_{\mu:\mu \unlhd \la}\chi_{\mu,\nu}^{(k)}\ d_{\mu\la}^{(k)}
\:.
$$ 
Note here that any partition $\mu$ dominated by a $k$-bounded partition~$\la$
is also $k$-bounded. 
As $D^{(k)}=(d_{\la\mu}^{(k)})_{\la,\mu\in P^{(k)}}$ and
$D^{(k)}_{(\al)}=(d_{\la\mu}^{(k)})_{\la,\mu<\al}$
are upper unitriangular integral matrices with
$$X^{(k)}=(D^{(k)})^t \X^{(k)}\:, \: X_{(\al)}=(D^{(k)}_{(\al)})^t \X^{(k)}_{(\al)}\;,$$
the claim is proved.

(2)
We have computed the determinant for $X_{(\al)}$
explicitly in Theorem~\ref{thm:small-basicsets}
to be
$$
\det X_{(\al)} = \prod_{\la < \al} a_{\la} \:,
$$
hence by (1) we get the formula for $\det  \X^{(k)}_{(\al)}$.

Note that $P^{(k)}=\{\al \in P(n) \mid \al < (1^{n-k-1},k+1)\}$,
and thus in our earlier notation  $X^{(k)}=X_{(1^{n-k-1},k+1)}$, 
giving also the assertion on $\det \X^{(k)}$
by Theorem~\ref{thm:small-basicsets}.
\qed

\medskip
\begin{rem}
{\rm 
As pointed out by a referee,
the formula for the determinant of the transition matrix
$\X^{(k)}$ between the $k$-Schur functions and the power sum
functions may also be obtained more directly.
For this, observe that the transition matrix between the $k$-Schur functions and the complete
symmetric functions is unitriangular \cite[eq. (6) and (7)]{LM},
and that the transition matrix between the complete
symmetric functions  associated
to $k$-bounded partitions and the $k$-bounded
power sums is triangular with the numbers $a_\la$, $\la\in P^{(k)}$,
on the diagonal.
}
\end{rem}

\medskip

For the dual coefficient matrix we deduce the following.
\begin{cor}
For $k\in \N$, we have
$$
\det \tilde\X^{(k)} = \prod_{\la \in P^{(k)}} b_{\la} \:.
$$
\end{cor}
\proof  For the matrices~$\tilde\X^{(k)}$, Theorem~\ref{thm:k-Schur}(2)
and the duality relation between $\X^{(k)}$ and  $\tilde\X^{(k)}$
immediately yield the assertion.
\qed

\medskip

Indeed, one observes an even better behavior in the data, analogous to the
phenomenon proved for the matrices $X^{(\al)}$ in Theorem~\ref{thm:smith},
namely, the Smith normal form satisfies
$$\S(\tilde\X^{(k)})=
\S(b_{\la}; \la \in P^{(k)})\:.$$

Dually to the refined observations on the (upper) principal
submatrices of $\X^{(k)}$,
for a given $\al \in P^{(k)}$
we also consider the (lower) principal submatrix
$$(\tilde\X^{(k)})^{(\al)}=
(\tilde\chi_{\la,\nu}^{(k)})_{\la,\nu \in P^{(k)}\atop \la,\nu\geq \al}$$
of $\tilde\X^{(k)}$.
This behaves like the corresponding submatrix in the character table,
again refining the observation made above, i.e., explicitly,
$$\S((\tilde\X^{(k)})^{(\al)})=
\S(b_{\la}; \la \in P^{(k)}, \la \geq \al)\:.$$

These final observations still need to be explained.

\bigskip

{\bf Acknowledgements.}
The authors would like to thank the Danish Research Council (FNU)
for the support of their collaboration.
We also express our thanks to Anne Schilling for sharing
her data on $k$-Schur functions at the meeting on
{\em Algebraic Combinatorixx} at the Banff International Research Station
(BIRS), 2011, and for pointing out a result by Lapointe and Morse.
Thanks go also to the referees for their helpful
suggestions and comments.


\end{document}